\documentclass{amsart}
\usepackage{amsmath}
\usepackage{amsfonts}

\newcommand{\dsum}{\displaystyle\sum}
\newcommand{\dprod}{\displaystyle\prod}

\DeclareMathOperator{\CEPL}{L}
\DeclareMathOperator{\tr}{tr}
\newcommand{\ext}[2]{#1/#2}
\newcommand{\Gal}[2]{\mathop\mathrm{Gal}(\ext{#1}{#2})}
\newcommand{\kro}[2]{\left( \frac{#1}{#2} \right)}

\newcommand{\Q}{{\bf Q}}

{
\theoremstyle{definition}
\newtheorem{definition}{Definition}
}

\newtheorem{theorem}{Theorem}

\begin{document}

\title[Quadratic irregularity of prime
numbers]{Comparison of algorithms to calculate quadratic irregularity of prime
numbers}

\author{Joshua Holden}
\address{Department of Mathematics and Statistics\\
University of Massachusetts at Amherst\\
Amherst, MA 01003, USA}
\curraddr{    Joshua Holden    \\    
    Mathematics Department\\
Duke University\\
Box 90320\\
Durham, NC 27708-0320, USA}

\email{holden@math.duke.edu}
\urladdr{http://www.math.duke.edu/\~{}holden}



\thanks{To appear in \emph{Mathematics of Computation}.}  
  


\begin{abstract}
In previous work, the author has extended the concept of regular and
irregular primes to the setting of arbitrary totally real number
fields $k_{0}$, using the values of the zeta function $\zeta_{k_{0}}$
at negative integers as our ``higher Bernoulli numbers''.
In the case where $k_{0}$ is a real quadratic field, Siegel presented
two formulas for calculating these zeta-values: one using entirely
elementary methods and one which is derived from the theory of modular
forms.  (The author would like to thank Henri Cohen for suggesting an
analysis of the second formula.)  We briefly discuss several
algorithms based on these formulas and compare the running time
involved in using them to determine the index of $k_{0}$-irregularity
(more generally, ``quadratic irregularity'') of a prime number.
\end{abstract}

\maketitle

\section{Definitions}
Let $k_{0}$ be a totally real number field, and let $p$ be an odd
prime.  Let $k_{1} =
k_{0}(\zeta_{p})$, where $\zeta_{p^n}$ will denote a primitive
$p^n$-th root of unity.  Let $\Delta = \Gal{k_{1}}{k_{0}}$, and let
$\delta = \left|\Delta\right|$.  Let $p^{e}$ be the largest power
of $p$ such that $\zeta_{p^{e}} \in k_{0}(\zeta_{p})$.

\begin{definition}
Let $\zeta_{k_{0}}$ be the zeta function for $k_{0}$.  We
say that $p$ is \emph{$k_{0}$-regular} if $p$ is relatively
prime to $\zeta_{k_{0}}(1-2m)$ for all integers $m$ such that $2 \leq  2m
\leq \delta - 2$ and also $p$ is relatively prime to $p^{e}
\zeta_{k_{0}}(1 - \delta)$.  The number of such zeta-values that
are divisible by $p$ will be the \emph{index of $k_{0}$-irregularity}
of $p$.
\end{definition}

According to a well-known
theorem of Kummer, $p$ divides the order of the class group of
$\Q(\zeta_{p})$ if and only if $p$ divides the numerator of a
Bernoulli number $B_{m}$ for some even $m$ such that $2 \leq m \leq
p-3$.  Such primes are called irregular; the others are called
regular.
In the setting we have described above, the author proved in his
thesis~(\cite{dissert}, see also~\cite{Holden99}), building on work of
Greenberg and Kudo, that
under a certain technical condition Kummer's criterion can be extended
to give information about whether $p$ divides the class group of
$k_{0}(\zeta_{p})$. To be exact, let $k_{1}^{+}$ denote
the maximal real subfield of $k_{1}$, which is equal to
$k_{0}(\zeta_{p}+\zeta_{p}^{-1})$.  Let $h(k_{1})$ denote the class
number of $k_{1}$ and $h^{+}(k_{1})$ denote the class number of
$k_{1}^{+}$.  It is known that $h^{+}(k_{1}) \mid h(k_{1})$; we let
the relative class number $h^{-}(k_{1})$ be the quotient.

\begin{theorem}[Greenberg, Holden] \label{GH}
Assume that no prime of the field $k_{1}^{+}$ lying over $p$ splits
in $k_{1}$.  Then $p$ divides $h^{-}(k_{1})$ if and only if $p$
is not $k_{0}$-regular.
\end{theorem}

As an application, we note that one common way of constructing
public-key cryptographic systems is to utilize the problem of finding
a discrete logarithm in some abelian group.
In order to make sure that the discrete logarithm
problem is computationally hard, one needs to know something about the
structure of the group involved, e.g. that it is divisible by a large
prime.  Theorem~\ref{GH} shows that if $p$ is a large $k_{0}$-irregular prime
and the conditions of the theorem are met, then the class group of $k_{0}(\zeta_{p})$
may be especially suitable for cryptography.  (One should
see~\cite{BP97} for more on the use of class groups in cryptography.)

For the case we consider, $k_{0}$ will be a real quadratic field
$\Q(\sqrt{D})$, with $D$ a positive fundamental discriminant.  For
such a $k_{0}$, we will say that primes are $D$-regular or have given
index of $D$-irregularity, and we will let the zeta function
$\zeta_{k_{0}}$ be also denoted by $\zeta_{D}$.  (More generally, we
may refer to the concept as ``quadratic irregularity''.)  In this case
$\delta$ will be equal to $p-1$ unless $D = p$, in which case $\delta
= (p-1)/2$.  Also, $e$ is always equal to $1$ when $p$ does not divide
the order of $k_{0}$ over $\Q$, which is true in this case since $p$
is odd.  For the condition in Theorem~\ref{GH} that no prime of the
field $k_{1}^{+}$ lying over $p$ splits in $k_{1}$ to be satisfied it
is sufficient that $p$ should not divide $D$, and we should also note
that since $p$ does not divide the degree of $k_{0}=\Q(\sqrt{D})$ over
$\Q$, a theorem of Leopoldt shows that $p$ divides $h(k_{1})$ if and
only if $p$ divides $h^{-}(k_{1})$.

In general, we will consider three cost models for the time of
multiplication: first using naive multiplication ($O(tt')$), second
using Sch\"onhage-Strassen fast multiplication or a similar method
($O(t \lg(t')^{O(1)}$), and third using a model where multiplication
(or addition) takes constant time regardless of the size of the
factors ($O(1)$).  We do not expect constant time multiplication to
occur asymptotically in the real world, but it can provide useful
estimates in situations where the
size of the numbers involved is small compared to the word size of the
actual computer in question.  (In these running time bounds, $t$
is the number of bits in the larger multiplicand and $t'$ the number
of bits in the smaller.)

\section{First formula}

Siegel's first formula to compute \mbox{$\zeta_{D}(1-2m)$}
for $m \geq 1$ an integer is analogous to the formula $\zeta(1-2m) = -B_{2m}/(2m).$
Using elementary methods, Siegel showed that similarly
\begin{equation} \label{eq1}
\zeta_{D}(1-2m) = \frac{B_{2m}}{4m^{2}} D^{2m-1} \sum_{j=1}^{D}
\chi(j) B_{2m}(j/D).
\end{equation}
Here $\chi(j) = \kro{D}{j}$, the Kronecker symbol, and $B_{2m}(j/D)$
indicates the $2m$-th Bernoulli polynomial evaluated at the fraction
$j/D$.  The Bernoulli polynomial $B_{r}(x)$ can be computed from the
Bernoulli numbers as
\begin{equation*} 
B_{r}(x) = \sum_{s=0}^{r}\binom{r}{s} B_{r-s} x^{s}.
\end{equation*}

It is not difficult to estimate the sizes of the numbers involved.  We
will assume throughout that $B_{m}$, $1 \leq m \leq M$, are
precomputed over the common denominator of the final result, and
stored in this fashion each has size $O(m (\lg m + \lg D))$ bits for a
total table size of $O(M^{2}(\lg M + \lg D))$ bits.  (The
precomputation does not in fact add to the asymptotic running time.)
The rational numbers $B_{2m}(j/D)$ can also be stored in $O(m(\lg m +
\lg D))$ bits, as can the total.  See~\cite{Holden98} for more
details.

A first attempt at an
algorithm based on (\ref{eq1}) might compute
 $B_{0}(\alpha),
\ldots, B_{M}(\alpha)$ naively from the formula.  The time taken for
this would be dominated by the powerings.
For $\alpha = a/b$
some rational number, the total time with naive multiplication would be
$O(M^{4} (\lg M + \lg a + \lg b)^{O(1)}).$
Using fast multiplication instead of naive multiplication would
improve this to
$O(M^{3} (\lg M + \lg a + \lg b)^{O(1)}),$
while with constant time multiplication we need only time
$O(M^{2}\lg M)$
regardless of $a$ and $b$.

However we can do better than this, using a cross between Horner's
method of evaluating polynomials and an algorithm used by Brent to
calculate Bernoulli numbers, as was previously discussed by the author
in~\cite{Holden98}.  This method gives a total time of
$O(M^{3} (\lg M + \lg a + \lg b)^{O(1)})$
using either constant or fast multiplication, and time
$O(M^{2})$ using constant time multiplication.

Using either of these algorithms to compute \mbox{$\zeta_{D}(1-2m)$}, $2 \leq
2m \leq M$, is then relatively straightforward.  Note that the
Kronecker symbol $\chi(j)$ can be computed in time $O(\lg^{2} D)$.  The
slower version of the algorithm
has time
$O(M^{4}D (\lg M + \lg D)^{O(1)})$
with naive multiplication,
$O(M^{3}D (\lg M + \lg D)^{O(1)})$
with fast multiplication, and
$O(M^{2}D (\lg M + \lg D)^{O(1)})$
with constant time multiplication.
The faster version runs
in time
$O(M^{3}D (\lg D + \lg M)^{O(1)})$
with either naive or fast multiplication (the $O(1)$ factor is
different, of course)
and again in time
$O(M^{2}D (\lg M + \lg D)^{O(1)})$
with constant time multiplication.

\section{Second formula}
Siegel's second formula is, as I said, derived from the theory of
modular forms.  In general, for $k_{0}$ a totally real number field as
above, it says that
$$
\zeta_{k_{0}}(1-2m) = - 2^{n} c_{2mn}^{-1}\sum_{l=1}^{r} c_{2mn,l}
s_{l}^{k_{0}}(2m),
$$
where $n=[k_{0}:\Q]$, $c_{2mn}=c_{2mn,0}$ and $c_{2mn,l}$
are rational integers
depending only on $2mn$ and $l$ (given by explicit formulas which we
will discuss),
$$r=\begin{cases}
    \lfloor mn/6 \rfloor & \text{if $2mn \equiv 2$ modulo $12$} \\
    \lfloor mn/6 \rfloor + 1 & \text{otherwise},
\end{cases}
$$
and $s_{l}^{k_{0}}$ is a sum over norms of ideals in the ring of
integers of $k_{0}$, namely
$$s_{l}^{k_{0}}(2m) = \sum_{\nu\in(\mathfrak{d})^{-1},\ \nu \gg
0,\ \tr(\nu)=l} \sigma_{2m-1}((\nu)\mathfrak{d}),$$
where
$$\sigma_{2m-1}(\mathfrak{A})=\dsum_{\mathfrak{B}\mid \mathfrak{A}}
N(\mathfrak{B})^{2m-1}$$
is a generalization of the usual sum of powers function
and $\mathfrak{d}$ is the different of $k_{0}$.  In the quadratic case
this all becomes much easier:
\begin{equation} \label{eq2}
\zeta_{D}(1-2m) = - 4  c_{4m}^{-1}\sum_{l=1}^{r} c_{4m,l} s_{l}^{D}(2m), \qquad
r= \lfloor m/3 \rfloor + 1,
\end{equation}
$$s_{l}^{D}(2m) = \sum_{\nu\in(\sqrt{D})^{-1},\ \nu \gg
0,\ \tr(\nu)=l} \sigma_{2m-1}((\nu\sqrt{D})),$$
and $s_{l}^{D}(2m)$
can also be expressed in terms of a purely arithmetic function
$e_{2m-1}(n)$, as follows:
$$s_{l}^{k_{0}}(2m) = \sum_{j\mid l}\chi_{D}(j) j^{2m-1}
e_{2m-1}((l/j)^{2}D)$$
and
$$e_{2m-1}(n) = \sum_{\substack{x^{2}\equiv n \pmod{4}\\ \left| x \right| \leq
\sqrt{n}}} \sigma_{2m-1}\left(\frac{n-x^{2}}{4}\right)$$
where
$$\sigma_{2m-1}(n)=\dsum_{d\mid n}d^{2m-1}$$
is the usual sum-of-powers function.  (See~\cite{Siegel68},
\cite{Zagier}, \cite{Cohen75} and~\cite{Cohen76} for more detailed
descriptions of these formulas, and for their derivations.)  The
coefficients $c_{4m,l}$ are most easily expressed as the coefficients
of a certain power series, and can be computed as needed without
adding to the asymptotic running time. We will give explicit formulas
for the power series and discuss its computation in
Section~\ref{coefficients}.  It is not hard to prove that in this form
$c_{4m,l}$ is of size $O(m)$.  The running time for calculating the
function $e_{2m-1}(n)$ is complicated by the need for factoring; we
use here an estimate based on the elliptic curve factoring method (we
would expect something very similar with any of the other standard
subexponential methods) to get an expected running time involving the
function $\CEPL(x)=e^{\sqrt{\log x \log \log x}}$.  Given this, we get
an expected running time to compute $e_{2m-1}(n)$ of
$$O(\sqrt{n}
\CEPL(n)^{1+o(1)}+(2m-1)^{2}\sqrt{n} \lg^{2} n)$$
using naive multiplication.
If we now applied
(\ref{eq2}) as written to compute all \mbox{$\zeta_{D}(1-2m)$}, $2 \leq 2m \leq
M$, we would get a running time of
$$O(M^{3}
\sqrt{D} \CEPL(M)^{O(1)} \CEPL(D)^{O(1)} \lg M + M^{5}
\sqrt{D} \lg M (\lg M + \lg D)^{O(1)}),$$
again using naive multiplication.

However, it is more efficient to rearrange the terms of the formula
as follows:
\begin{eqnarray}
\nonumber
\zeta_{D}(1-2m) &= & -4 c_{4m}^{-1}\sum_{l=1}^{r} c_{4m,l} s_{l}^{k_{0}}(2m) \\
\nonumber
&=& -4 c_{4m}^{-1} \sum_{l=1}^{r} c_{4m,l} \sum_{j\mid l}\chi_{D}(j) j^{2m-1}
e_{2m-1}((l/j)^{2}D) \\
\label{eq3}
&=& -4 c_{4m}^{-1} \sum_{k=1}^{r} \left( \sum_{j=1}^{\lfloor r/k \rfloor}
\chi_{D}(j) j^{2m-1} c_{4m,jk} \right) e_{2m-1}(k^{2}D)
\end{eqnarray}
This rearrangement of the formula requires fewer calls to compute
$e_{2m-1}$ by a factor of $\lg m$.  Using this version of the formula,
the time necessary to compute all \mbox{$\zeta_{D}(1-2m)$}, $2 \leq 2m \leq
M$, using naive multiplication is
$$O(M^{3} \sqrt{D} \CEPL(M)^{O(1)} \CEPL(D)^{O(1)}+
M^{5} \sqrt{D} (\lg M + \lg D)^{O(1)}).$$
This is much worse than the best algorithm based on
(\ref{eq1}) in terms of $M$, but it is better in terms of
$D$.  Also, except for one final division by $c_{4m}^{-1}$, all of
the arithmetic in this formula deals only with rational integers;
unlike the previous formulas.  Note that the first term comes from the
factoring process, while the second term comes from
multiplications.

It should be noted that the asymptotic running time of this algorithm is greatly
improved by using Sch\"onhage-Strassen fast multiplication or constant
time multiplication, in which
cases the second term becomes smaller than the first and the running
time becomes
$$O(M^{3} \sqrt{D} \CEPL(M)^{O(1)} \CEPL(D)^{O(1)})$$
This is still worse than using (\ref{eq1}) in terms of $M$, but
only by a subexponential factor.

It should also be noted that (\ref{eq2}) and~(\ref{eq3}) also present
opportunities for time savings when computing zeta-values for multiple
$D$ in the same range of $M$, at a sacrifice of memory space.  The
controlling factor in the speed of the algorithm is the number of
times that $\sigma_{2m-1}(n)$ must be calculated.  Note that in computing
all \mbox{$\zeta_{d}(1-2m)$}, $5 \leq d \leq D$, there can only be
$O(m^{2}D)$ different values of $n$.  However, following the algorithm
strictly, we would ordinarily make $O(m^{2}D^{3/2})$ calls to the
subroutine that calculates this function.

Thus if we compute all \mbox{$\zeta_{d}(1-2m)$}, $5 \leq d \leq D$,
storing values of $\sigma_{2m-1}(n)$ as we compute them, and then repeat
this process for each $m$ in the range $2 \leq 2m \leq M$,  the running
time should be $O(D \CEPL(D)^{O(1)})$ in terms of $D$, rather than
$O(D^{3/2} \CEPL(D)^{O(1)})$ as one would obtain following the
algorithm strictly.
This compares very favorably with the time of $O(D^{2})$ in terms of
$D$ which holds for algorithms using (\ref{eq1}).

Since the exponent $2m-1$ used in the $\sigma_{2m-1}(n)$ function
changes as $m$ does, we can dispose of the table when we change $m$.
The table that we need to keep requires at most $O(M^{3}D(\lg M + \lg
D))$ bits of storage, which could be a significant barrier.  More
efficient storage of the important information may be
valuable here; we will discuss this somewhat more in
Section~\ref{conclusion}.

\section{Computing the numbers $c_{4m,l}$}
\label{coefficients}

The integers $c_{4m,l}$ are defined as follows.
Let
$$G_{k} =1-\dfrac{2k}{B_{k}}\dsum_{n=1}^{\infty}
\sigma_{k-1}(n) \, q^{n}$$
be the (normalized) Eisenstein series of order $k$ for $k=6, 10$, and
$14$.  (For the general $c_{2mn,l}$ one also needs $k=0, 4$, and $8$.)
Let
$$\begin{aligned}[t]
\Delta &= q \dprod_{n=1}^{\infty}(1-q^{n})^{24}\\
&= q \left(\dsum_{n=0}^{\infty} (-1)^{n} (2n+1) \,
q^{n(n+1)/2}\right)^{8}.
\end{aligned}$$
be the discriminant series.
Let $r= \lfloor m/3 \rfloor + 1$ as before, and let
$$T_{4m}=G_{12r-4m+2}\Delta^{-r} = \dsum_{n=-r}^{\infty} c_{4m, -n} \,
q^{n}.$$
Then $c_{4m}=c_{4m, 0}$, and the other $c_{4m,l}$ for $1 \leq l \leq
r$ can also be read off as coefficients of $T_{4m}$.  Luckily, the
expression $12r-4m+2$ only takes on the values $6$, $10$, and $14$.
(In the more general case we can define $T_{2mn}$ similarly; the
expression $12r-2mn+2$ can take on the values $0$, $4$, and $8$ in
addition to those above.)

The best algorithm known to the author for calculating these
coefficients goes roughly as follows.  At the start of the
computations for \mbox{$\zeta_{D}(1-2m)$}, $2 \leq 2m \leq M$, calculate
$G_{6}$, $G_{10}$, $G_{14}$, and $\Delta^{-1}$ with the maximum number
of coefficients necessary (about $M/12$).  Instead of trying to
compute all of the needed series $\Delta^{-r}$ at once, we calculate
it as a running product which only needs to be updated when $r$
changes.  Then, whenever $m$ changes, we multiply truncated versions
$G_{12r-4m+2}$ and $\Delta^{-r}$ (with about $m/6$ coefficients each)
to find the required coefficients of $T_{4m}$.

The series $\Delta^{-1}$ can also be expressed as
$$\begin{aligned}[t]
\Delta^{-1} &= q \prod_{n=1}^{\infty}(1-q^{n})^{-24}\\
&= q \left(\prod_{n=1}^{\infty}\frac{1}{(1-q^{n})}\right)^{24}\\
&= q \left(\sum_{n=0}^{\infty} p(n)\,q^{n}\right)^{24}
\end{aligned}$$
where $p(n)$ takes on integer values and is well-known as the
partition function from additive number theory.  Hardy and Ramanujan
proved an asymptotic expression for $p(n)$ which shows that $\lg p(n)$
is of order $\sqrt{n}$.  (See, for example, Chapter~14 of~\cite{Apostol}.)  From
this it is easy to show that the coefficients of $\Delta^{-r}$ take at
most $O(M^{1.5})$ bits of storage each, as do the coefficients of
$G_{12r-4m+2}$ and $T_{4m}$.  Thus
the storage required for all of the computation of $c_{4m, l}$
necessary for a fixed $m$ is $O(M^{2.5})$.  The resulting table, which
is of course independent of $d$, can be
used to compute \mbox{$\zeta_{d}(1-2m)$} for any range of $d$ and then
disposed of when $m$ is changed.

As far as the time for this algorithm is concerned, computing the
tables necessary for all \mbox{$\zeta_{d}(1-2m)$}, $2 \leq 2m \leq M$,
$5 \leq d \leq D$, takes time $O(M^{5})$ with naive multiplication of
integers and also of polynomials.  By using FFT methods to multiply
polynomials the time can be reduced to $O(M^{3.5} \lg M \lg \lg M)$
with fast multiplication of integers and $O(M^{2} \lg M)$ with
constant time multiplication of integers.  (The use of FFT methods
here was suggested by A.O.L. Atkin and Will Galway.)  This is within
our previously established time bounds in the naive and constant time
cases; however in the fast multiplication case it could add to the
total asymptotic time in terms of $M$, which for the previous parts of
the algorithm was established as $O(M^{3} \CEPL(M)^{O(1)})$ in terms
of $M$.  On the other hand it should be noted that these calculations
only need to be done once per value of $m$ no matter how many values
of $d$ one is examining.  Also the constant involved in the $O(M^{3.5}
\lg M \lg \lg M)$ seems to be quite good in practice compared to that
in the $O(M^{3} \CEPL(M)^{O(1)})$.

\section{Summary and future work}

\label{conclusion}

Table~\ref{algstable} presents the various algorithms, for comparison.
We present the asymptotic order of the running time to compute
\mbox{$\zeta_{D}(1-2m)$}, $2 \leq 2m \leq M$, for a given $D$, using
the three methods of multiplication discussed earlier.


\begin{table}

    \caption{Comparison of algorithms for \mbox{$\zeta_{D}(1-2m)$}, $2 \leq 2m \leq M$}
    \label{algstable}
\begin{tabular}{lll}
Equation &&\\
used & Multiplication & Time order \\
\hline
(\ref{eq1})  & naive & $M^{4}D (\lg M + \lg D)^{O(1)}$\\
(\ref{eq1})  & fast & $M^{3}D (\lg M + \lg D)^{O(1)}$\\
(\ref{eq1})  & constant &  $M^{2} D (\lg M + \lg D)^{O(1)}$\\

(\ref{eq1}) from~\cite{Holden98} & naive & $M^{3}D (\lg D + \lg M)^{O(1)}$ \\
(\ref{eq1}) from~\cite{Holden98} & fast & $M^{3}D (\lg D + \lg M)^{O(1)}$ \\
(\ref{eq1}) from~\cite{Holden98} & constant & $M^{2} D (\lg M + \lg D)^{O(1)}$ \\

(\ref{eq2}) & naive & $M^{3} \sqrt{D} \CEPL(M)^{O(1)} \CEPL(D)^{O(1)} \lg
M$ \mbox{\qquad} \\
& & \hfill  $+ M^{5} \sqrt{D} (\lg M + \lg D)^{O(1)}$ \\
(\ref{eq2}) & fast & $M^{3} \sqrt{D} \CEPL(M)^{O(1)}
    \CEPL(D)^{O(1)} \lg M$ \\
(\ref{eq2}) & constant & $M^{3} \sqrt{D} \CEPL(M)^{O(1)}
    \CEPL(D)^{O(1)} \lg M$ \\

(\ref{eq3}) & naive & $M^{3} \sqrt{D} \CEPL(M)^{O(1)} \CEPL(D)^{O(1)}$
\mbox{\qquad}\\
& & \hfill $+M^{5} \sqrt{D} (\lg M + \lg D)^{O(1)}$ \\
(\ref{eq3}) & fast & $M^{3} \sqrt{D} \CEPL(M)^{O(1)}
    \CEPL(D)^{O(1)}$ \\
(\ref{eq3}) & constant & $M^{3} \sqrt{D} \CEPL(M)^{O(1)}
    \CEPL(D)^{O(1)}$\\
\end{tabular}
\end{table}

The factor of $\lg M$ in the times for algorithms based on~(\ref{eq2}) has been
included to emphasize that these algorithms are slower than those based
on~(\ref{eq3}), even though the factor of $\lg M$ could be
absorbed into that of $\CEPL(M)^{O(1)}$.

We also provide, in Tables~\ref{timetable} and~\ref{timetable2},
tables of actual timings for some of the algorithms, using naive
multiplication.  The times were measured on a Sun SPARC Ultra-1
computer using the GP-Pari interpreted language.
(See~\cite{PARImanual}.)  The data computed by these programs will be
analyzed in a future paper.

Table~\ref{timetable} measures the time to compute
\mbox{$\zeta_{D}(1-2m)$}, $2 \leq 2m \leq M$, for a given $D$, in
hours, minutes, and seconds.  The number in parentheses indicates
the size of the stack used, in bytes.  These numbers are only a very
rough guide to the actual amount of memory used.

\begin{table}

\caption{Calculating \mbox{$\zeta_{D}(1-2m)$}, $2 \leq 2m \leq M$}
    \label{timetable}

\begin{tabular}{|r|r|rr|rr|}
\hline
$D$ & $M$  & time (\ref{eq1}) from~\cite{Holden98} &(stacksize)
& time (\ref{eq3}) &(stacksize)\\
\hline
5   &  100 & 3.155    & (10M)    &  .838 &(10M)\\
101 & 100  & 55.561   & (10M)    & 3.758 & (10M)\\
501 & 100  & 4:50.282  & (10M)   & 12.480 &(10M)\\
1001 & 100 & 10:52.411  & (10M)  & 20.670 & (10M)\\
5001 & 100 & 48:27.107  &(12M)   & 1:43.548 &(10M)\\
\hline
5 & 500    & 2:05.670  & (4M)    & 8:11.603  &(4M) \\
101 & 500  & 42:38.612  &(4M)   & 1:18:26.615 &(4M) \\
501 & 500   & 3:48:18.615& (8M)  & 4:45:20.438 &(4M) \\
1001 & 500  & 7:49:56.048  &(12M) & 7:49:00.783 &(4M) \\
\hline
5 & 1000   & 10:05.903 & (4M)    & 3:55:46.908 &(4M) \\
\hline
5 & 2000   & 1:14:01.992 & (16M)  & 118:17:46.020 &(4M) \\

\hline
\end{tabular}
\end{table}

Table~\ref{timetable2} measures the time to compute
\mbox{$\zeta_{d}(1-2m)$}, $2 \leq 2m \leq M$, $5 \leq d \leq D$.  We
use the algorithm based on~(\ref{eq3}), both with and without keeping
a table of $\sigma_{2m-1}(n)$ as described earlier.  The units and
stack size numbers should be interpreted as in Table~\ref{timetable}.

\begin{table}

\caption{Calculating \mbox{$\zeta_{d}(1-2m)$}, $2 \leq 2m
\leq M$, $5 \leq d \leq D$}
    \label{timetable2}

\begin{tabular}{|r|r|rr|rr|}
\hline
$D$ & $M$  & time (\ref{eq3}) & (stacksize)
& time (\ref{eq3}) with table & (stacksize)\\
\hline
100 & 100  & 1:19.638  &  (4M)    & 1:27.983 & (4M)\\
500 & 100  & 18:34.377  &(4M)   & 13:57.166 &(8M)\\
1000 & 100 &  1:03:50.464 &(4M)  &  37:52.552 &(12M)\\
5000 & 100 &  26:39:46.109& (4M)   &  7:10:44.870 &(64M)\\

\hline
\end{tabular}
\end{table}

As one can see, memory usage for algorithms based on~(\ref{eq3}) with
a table of values $\sigma_{2m-1}(n)$ goes up quite quickly, and even
so a lot of redundant work is being done.  For one thing, the same
numbers $n$ will have to be factored repeatedly for different values
of $2m-1$, but they will lead to different values of
$\sigma_{2m-1}(n)$ so the actual factorization would need to be stored
and not just a function value.  This would be even more memory
intensive.  Carl Pomerance has suggested some approaches to this,
including using a cache rather than a complete table, and storing only
the largest prime factor rather than a complete factorization.  The
question of a cache leads naturally to the question of which numbers
will appear as integer values of $(k^{D}-x^{2})/4$ as $k$ and $D$
vary, and how often.  Henri Cohen, in~\cite{Cohen76}, gives some
variations on Siegel's formula which could cut down on the number of
times $\sigma_{2m-1}(n)$ needs to be computed.  Algorithms based on
these might provide a linear speedup over the algorithms presented
here, but the asymptotic behavior would probably be the same.

The anonymous reviewer has suggested that some of the
arithmetic could possibly be speeded up by the use of modular
techniques and the Chinese remainder theorem.  This certainly deserves
more consideration.

Another prospect for future work is the analysis of ``first-hit''
versions of these algorithms, namely determining how long we should
expect to search before finding, say, the first $D$-irregular prime
larger than a certain bound for $D$ in a given range.  This might be
particularly useful for cryptographic applications, in which we would
be explicitly looking for a class group (or a small number of them) with
a hard discrete logarithm problem.  This will be explained in more detail
in~\cite{mathcomp2}.

\section*{Acknowledgements}
The author would like to thank A.O.L. Atkin, Will Galway, and the
members of the NMBRTHY electronic mailing list for their very helpful
suggestions, and Robert Harley for generating some of the precomputed
tables used in the computations reported in this paper.  He would also
like to especially thank Johannes Buchmann, Henri Cohen, and Carl
Pomerance for suggestions and encouragement, and Gary Walsh for his
encouraging remarks on an early version of this paper.  Finally, he
would like to thank the anonymous reviewer for encouragement and
several helpful suggestions.

\newcommand{\SortNoop}[1]{}
\providecommand{\bysame}{\leavevmode\hbox to3em{\hrulefill}\thinspace}

\end{document}